\documentclass[lefttitle,nopreprintline,11pt]{elsarticle}

\usepackage[T1]{fontenc}
\usepackage{amsmath,amssymb,amsthm}
\usepackage[hidelinks]{hyperref}

\newcommand{\dd}{\,\mathrm d}
\newcommand{\e}{\mathrm e}
\newcommand{\ii}{\mathrm i}

\journal{Advances in Mathematics}

\def\dochead#1{\csgappto{@title}{\scalebox{0.8}{#1}\par}}
\patchcmd{\pprintMaketitle}{\hrule\vskip12pt}{\vspace*{-24pt}}{}{}
\patchcmd{\MaketitleBox}{\hrule\vskip12pt}{\vspace*{-24pt}}{}{}
\DeclareDocumentCommand{\originaldoi}{m}{\csgappto{@cornotes}{\insert\footins{\footnotesize#1}}}

\newtheorem{theorem}{Theorem}
\newtheorem{lemma}[theorem]{Lemma}

\theoremstyle{remark}

\DeclareMathOperator{\Mp}{Mp}
\newcommand{\R}{\mathbb{R}}
\newcommand{\cS}{\mathcal{S}}
\newcommand{\eps}{\varepsilon}
\newcommand{\abs}[1]{\lvert#1\rvert}
\newcommand{\ip}[2]{\langle #1,#2\rangle}

\begin{document}

\begin{frontmatter}
    \dochead{Corrigendum}
    \title{Corrigendum to "Sparse Gabor representations of metaplectic operators: controlled exponential decay and Schr\"odinger confinement" [Adv. Math. 490 (2026) 110828]}
    \originaldoi{DOI of original article: \href{https://doi.org/10.1016/j.aim.2026.110828}{10.1016/j.aim.2026.110828}.}
    \author[unito]{Elena Cordero}
    \ead{elena.cordero@unito.it}
    \author[usi]{Gianluca Giacchi}
    \ead{gianluca.giacchi@usi.ch}
    \author[unito]{Edoardo Pucci}
    \ead{edoardo.pucci@unito.it}
    \author[polito]{S. Ivan Trapasso\corref{cor1}}
    \ead{salvatoreivan.trapasso@polito.it}
    \cortext[cor1]{Corresponding author.}
    \affiliation[unito]{organization={Dipartimento di Matematica, Universit\`a di Torino}, addressline={via Carlo Alberto 10}, city={Torino}, postcode={10123}, country={Italy}}
    \affiliation[usi]{organization={Faculty of Informatics, Universit\`a della Svizzera italiana}, addressline={Via la Santa 1}, city={Lugano}, postcode={6962}, country={Switzerland}}
    \affiliation[polito]{organization={Dipartimento di Scienze Matematiche "G. L. Lagrange", Politecnico di Torino}, addressline={corso Duca degli Abruzzi 24}, city={Torino}, postcode={10129}, country={Italy}}
\end{frontmatter}

The purpose of this note is to correct some quantitative estimates and computations in \cite{CGPT}. The proofs of Theorem~1.3 and Lemma~3.6 there contain gaps, and consequently some of the numerical decay coefficients that depend on them need adjustment. Some identities concerning Hamiltonian flows and Euler decompositions in Section~6 also require correction. Theorem~1.2 and the qualitative conclusions of the general metaplectic estimates in \cite{CGPT}, namely Gabor sparsity and exponential confinement for quadratic Schr\"odinger evolutions, are not affected. We use the notation of \cite{CGPT}, write $p=1/s$, and retain its Euler matrices $D',D''$.

\paragraph{Theorem 1.3}
Two steps of the proof require more care: the radial weight is split into one-variable factors without accounting for the loss in its coefficient that this entails, and the last step divides by an upper bound for the window, which does not control the window from below. In fact, the coefficient stated in the conclusion cannot be retained in general: for $s=1/2$ and $0<\eps<\pi$, the functions $g(t)=\e^{-\pi t^2}\in\cS^{1/2}_{1/2,\eps}(\R)$ and $f(t)=\e^{-at^2}$, where $a=\pi\eps/(2\pi-\eps)<\eps$, satisfy
\[ \abs{V_gf(x,\xi)}=\left(\frac{\pi}{a+\pi}\right)^{1/2}\e^{-\eps x^2/2-(\pi-\eps/2)\xi^2}\leq\e^{-\eps\abs{(x,\xi)}^2/2}, \]
while $f\notin\cS^{1/2}_{1/2,\eps}(\R)$. The following result is the appropriate replacement; note that it does not require fixed-rate bounds on the window.

\begin{theorem}\label{thm:13}
Let $s\geq1/2$, $g\in\cS(\R^d)\setminus\{0\}$ and $f\in\cS'(\R^d)$. If $\abs{V_gf(z)}\leq C\e^{-A\abs z^{1/s}}$ for some $A,C>0$ and all $z\in\R^{2d}$, then $f\in\cS^s_{s,\rho}(\R^d)$ for every $0<\rho<A$. At $s=1/2$ one also has $f\in\cS^{1/2}_{1/2,A}(\R^d)$.
\end{theorem}

\begin{proof}
Fourier inversion in the second STFT variable yields, for $g(y_0)\neq0$, 
\[ \abs{f(t)}\abs{g(y_0)}\leq\int_{\R^d}\abs{V_gf(t-y_0,\xi)}\dd\xi. \]
Fix now $\rho<\theta A<A$; since $\theta\abs x^p+(1-\theta)\abs\xi^p\leq\abs{(x,\xi)}^p$, integration gives $\abs{f(t)}\lesssim\e^{-\rho\abs t^p}$, and the Fourier-covariance identity provides the same bound for $\widehat f$, hence smoothness. For $p=2$ and $f\neq0$ we use the STFT symmetry $\abs{V_fg(z)}=\abs{V_gf(-z)}$ and the first assertion to infer that $g$ and $\widehat g$ extend to nonzero entire functions. Choose $y_1,\ldots,y_{d+1}$ with $g(y_j)\neq0$ and $0$ in their convex hull; for every $t$ some $j$ satisfies $t\cdot y_j\leq0$, so that $\abs{t-y_j}^2\geq\abs t^2$ and
\[ \abs{f(t)}\leq \frac{C}{\min_j\abs{g(y_j)}}\left(\frac{\pi}{A}\right)^{d/2}\e^{-A\abs t^2}. \]
The same argument with $\widehat g$ yields the endpoint bound for $\widehat f$.
\end{proof}

For $s>1$ the hypothesis of \cite[Theorem~1.3]{CGPT} has $A=2^{1-1/s}\eps>\eps$, so that choosing $\rho=\eps$ recovers its conclusion as stated. At $s=1$ one has $A=\eps$, and Theorem~\ref{thm:13} gives the conclusion for every coefficient $\rho<\eps$.

\paragraph{Lemma 3.6}
The coefficient $\psi_s(\eps)$ in \cite[Lemma~3.6]{CGPT} is too large in general: at $s=1/2$, $d=1$ and $\Sigma=I$ the convolution equals $(\pi/(2\eps))\exp(-\eps\abs v^2/2)$, which decays with coefficient $\eps/2$ rather than $2\eps$. The following statement replaces that lemma.

\begin{lemma}\label{lem:36}
Set $\kappa_s=2^{1-1/s}$ for $1/2\leq s\leq1$ and $\kappa_s=1$ for $s\geq1$. For every $\eps>0$ and $0<\eta<\kappa_s\eps$,
\[ \int_{\R^{2d}}\e^{-\eps\abs{D'u-v}^{1/s}}\e^{-\eps\abs{D''u}^{1/s}}\dd u\leq C_{s,d,\eps,\eta}\e^{-\eta\abs{D''v}^{1/s}},\qquad v\in\R^{2d}, \]
with a constant independent of $\Sigma=\operatorname{diag}(\sigma_1,\ldots,\sigma_d)$, $\sigma_j\geq1$. At $s=1/2$ one may take $\eta=\eps/2$ and $C=(\pi/\eps)^d$.
\end{lemma}

\begin{proof}
Since $D''$ is a contraction and $D''D'=\Sigma^{-1}\oplus\Sigma^{-1}$, we have $\abs{D''v}\leq\abs{D'u-v}+\abs{D''u}$, and the inequality $a^p+b^p\geq\kappa_s(a+b)^p$ gives the claimed exponential with $\eta=\theta\kappa_s\eps$, $0<\theta<1$. The remaining factor is bounded by $\exp(-(1-\theta)\eps(\abs{v_2-u_2}^p+\abs{u_1}^p))$, whose integral is uniform in $\Sigma,v$. If $p=2$, completing the square with $a_j=\sigma_j^{-1}\leq1$ gives the exact integral
\[ \left(\frac{\pi}{\eps}\right)^d\prod_{j=1}^d\frac{\exp[-\eps(v_j^2+a_j^2v_{d+j}^2)/(1+a_j^2)]}{1+a_j^2}, \]
which is bounded by $(\pi/\eps)^d\exp(-\eps\abs{D''v}^2/2)$.
\end{proof}

\paragraph{Consequences}
In the regime $s\geq1$, the conclusions of Theorem~1.4, Corollaries~1.5 and 1.6, Proposition~5.5 and Corollary~5.6 of \cite{CGPT} remain valid as stated: the coefficients in those statements lie strictly within the bounds of Lemma~\ref{lem:36}, and the open parameter ranges allow the use of the strict converse in Theorem~\ref{thm:13}. For $1/2\leq s<1$, the corresponding numerical conclusions are to be replaced as follows.

\begin{theorem}\label{thm:14}
Let $1/2\leq s<1$, $\eps>0$, and $g,\gamma\in\cS^s_{s,\eps}(\R^d)$. Put $\delta=\eps/4$ if $s=1/2$; if $1/2<s<1$, let $0<\delta<2^{1-3/(2s)}\eps$. There is $C>0$ such that, for every $\widehat S\in\Mp(d,\R)$ and every Euler decomposition of its projection $S$,
\[ \abs{\ip{\widehat S\pi(z)g}{\pi(w)\gamma}}\leq C\det(\Sigma)^{-1/2}\e^{-\delta\abs{D'U(w-Sz)}^{1/s}},\qquad z,w\in\R^{2d}. \]
Consequently, for $f,g\in\cS^s_{s,\eps}(\R^d)$ and $\widehat S\in\Mp(d,\R)$,
\[ \widehat Sf\in\cS^s_{s,\rho}(\R^d),\qquad V_g(\widehat Sf)\in\cS^s_{s,\rho}(\R^{2d}), \]
with $\rho=\eps\sigma_{\min}^2/4$ if $s=1/2$, and with every $0<\rho<2^{1-3/(2s)}\eps\sigma_{\min}^{1/s}$ if $1/2<s<1$.
\end{theorem}

\begin{proof}
By \cite[Theorem~1.2 and (21)]{CGPT}, the Wigner distributions of $g,\gamma$ decay with coefficient $2^{1/s}\delta_0$, where $\delta_0=\eps/4$ for $s=1/2$ and $\delta<\delta_0<2^{1-3/(2s)}\eps$ otherwise. The Wigner reduction in \cite[Section~4.3]{CGPT} (which remains valid once Lemma~3.6 is replaced by Lemma~\ref{lem:36}) then gives the stated Gabor bound because $\kappa_s2^{1/s}=2$, while in the Gaussian case we use the endpoint of that lemma. Setting $z=0$ and using $\abs{D'Uw}\geq\sigma_{\min}\abs w$ yields an STFT bound with coefficient $\delta\sigma_{\min}^{1/s}$. Theorem~\ref{thm:13} then gives the estimate for $\widehat Sf$, including the Gaussian endpoint. The Fourier transform of $V_g(\widehat Sf)$ has modulus $\abs{\widehat Sf(-\eta)}\abs{\widehat g(\omega)}$, and $\rho\leq\eps$ gives the remaining membership.
\end{proof}

Accordingly, for $1/2\leq s<1$ Theorem~\ref{thm:14} supersedes the coefficients in Theorem~1.4, Corollaries~1.5 and 1.6, Propositions~5.1 and 5.8, Corollary~5.3 and the related remarks in \cite{CGPT}. We point out that for $s=1/2$, $d=1$ and $t=0$ the value $\eps(0)=2\eps$ in Corollary~1.6 would place $\e^{-\eps x^2}$ in $\cS^{1/2}_{1/2,2\eps}(\R)$, which is not the case for $0<\eps\leq\pi$; the value $\eps/(4d)$ mentioned after that corollary, which differs from the one given by its displayed formula, is superseded as well. Finally, the comparison with \cite{RR,Vemuri} should be read with an additional factor $\pi$ in our normalization: an initial rate $\pi\tanh(2\alpha)$ yields the all-time one-dimensional harmonic-oscillator rate $\pi\tanh\alpha$ in \cite{RR}, and every smaller rate in \cite{Vemuri}.

\paragraph{Section 5.3}
In \cite[Section~5.3]{CGPT} it is stated that (23) implies (22) with $a'=a/2$ and $b'=b/2$. This is not the case, and in fact no implication of this form can hold: the modulus of $V_ff$ is invariant under time-frequency shifts, since $\abs{V_{\pi(w)f}(\pi(w)f)}=\abs{V_ff}$ for every $w\in\R^{2d}$, whereas the bounds in (22) are not. For instance, $f(t)=\e^{-\pi(t-1)^2}$ satisfies
\[ \abs{V_ff(x,\xi)}=2^{-1/2}\e^{-\pi(x^2+\xi^2)/2}, \]
but $\e^{\pi t^2}\abs{f(t)}=\e^{2\pi t-\pi}$ is unbounded. The marginal argument given there actually applies to the Wigner distribution $Wf$, which by (13) is a rescaling of $V_{\mathcal If}f$ rather than of $V_ff$. Accordingly, the correct statement reads as follows: if $f\in L^2(\R^d)$ satisfies
\[ \abs{V_{\mathcal If}f(x,\xi)}\lesssim\e^{-a'\abs x^2-b'\abs\xi^2},\qquad x,\xi\in\R^d, \]
then for $s=1/2$ $f$ satisfies (22) with $a=2a'$ and $b=2b'$. Proposition~5.7 is not affected.

In the chirped-Gaussian example that follows, the Fourier transform of $f(x)=\e^{\ii\pi\abs x^2-a\abs x^2}$ is
\[ \widehat f(\xi)=\left(\frac a\pi-\ii\right)^{-d/2}\e^{-\pi^2\abs\xi^2/(a-\ii\pi)}, \]
so that (22) holds with $b=a\pi^2/(a^2+\pi^2)$ rather than $b=\pi^2/(a^2+\pi^2)$. The conclusion of the example is unchanged. Along the diagonal one has $\abs{V_ff(x,x)}=(2a/\pi)^{-d/2}\e^{-a\abs x^2/2}$, which is not bounded by a multiple of $\e^{-(a+b)\abs x^2/2}$ for any $b>0$; hence (22) does not imply (23) with $a'=a/2$ and $b'=b/2$.

\paragraph{Section 6}
Some explicit computations in \cite[Section~6]{CGPT} require the following corrections; none of them affects the results of Sections~1 to 5. In the first scenario of Section~6.2 the case $\omega=0$ must be excluded, since the flow is then the shear $\begin{pmatrix}1&t/m\\0&1\end{pmatrix}$, that is, the free particle of Section~6.1 with $t/m$ in place of $2t$; accordingly, in the higher-dimensional oscillator each block with $\omega_j=0$ is a free shear. For $m\omega=\pm1$ the flow is
\[ S_t=\begin{pmatrix}\cos(\omega t)&\operatorname{sgn}(m\omega)\sin(\omega t)\\-\operatorname{sgn}(m\omega)\sin(\omega t)&\cos(\omega t)\end{pmatrix}, \]
and an Euler decomposition is given by $U_t=D_t=I$ and $V_t=S_t$. In the second scenario the condition should read $\omega t=\pi/2+k\pi$, $k\in\mathbb{Z}$, which is the case in which $\sin(\omega t)=(-1)^k$; the subsequent computation is correct. In Section~6.3 the coefficient $\omega$ of the magnetic term in (30) should be replaced by $\omega/2$, consistently with the flow $S_t$ displayed there, and the factor $\frac1{m\omega}-\omega$ in the third scenario should read $\frac1{m\omega}-m\omega$. In Section~6.4.2 the middle factor in the decomposition of $\mathcal D_E$ should be $D^{-1}\oplus D$, and an Euler decomposition is obtained from it by exchanging $x_j$ and $\xi_j$ through the symplectic rotation of the $(x_j,\xi_j)$ plane whenever the $j$-th singular value of $E$ exceeds $1$, followed by a reordering. In Section~6.4.3 the normalizing factors should be $\rho_\pm=(\mu^2+(\lambda_\pm-1)^2)^{1/2}$, the norms of the eigenvectors $(\mu,\lambda_\pm-1)$ of $V_\mu V_\mu^\top$; with the values stated there the matrix $U$ is not orthogonal.

\paragraph{Declaration of generative AI and AI-assisted technologies in the writing process}
The authors used AI tools (ChatGPT, Claude) to assist with detecting inaccuracies in the original article and drafting this corrigendum. They reviewed and edited the content as needed and take full responsibility for the content of this manuscript.

\end{document}